\documentclass[twoside]{article}
\usepackage{amssymb,amsfonts}
\usepackage{graphicx}

\makeatletter
\def\ps@myheadings{\let\@mkboth\@gobbletwo
  \def\@oddhead{\parbox{\textwidth}%
     {\mbox{}\hfill\rm\thepage\hfill\mbox{}}}%
  \def\@oddfoot{}%
  \def\@evenhead{\parbox{\textwidth}%
    {\mbox{}\hfill\rm\thepage\hfill\mbox{}}}%
  \def\@evenfoot{}%
 \def\chaptermark##1{}\def\sectionmark##1{}\def\subsectionmark##1{}}
\makeatother

\newtheorem{theorem}{Theorem}[section]        
\newtheorem{definition}[theorem]{Definition}     
\newtheorem{lemma}[theorem]{Lemma}             
\newtheorem{remark}[theorem]{Remark}           

\newtheorem{proposition}[theorem]{Proposition}
\newenvironment{proof}{{\it Proof.}}{\mbox{\ }\hfill 
$\square$\vskip3mm}         
\newenvironment{proofof}[1]{{\it #1.}}{\mbox{\ }\hfill 
$\square$\vskip3mm}         
\newcommand{\R}{{\mathbb R}}
\newcommand{\C}{{\mathbb C}}

\renewcommand{\epsilon}{\varepsilon}
\newcommand{\N}{{\mathbb N}}

\newcommand{\Cal}[1]{{\cal #1}}
\newcommand{\eqnref}[1]{{\rm (\ref{#1})}}
\newcommand{\mysection}[2]{\section*{{\large\bf #1. #2}}%
\setcounter{equation}{0}\setcounter{section}{#1}%
\setcounter{theorem}{0}}

\newcommand{\myappendix}[1]{\section*{{\large\bf Appendix. #1}}
\setcounter{equation}{0}
\setcounter{theorem}{0}
\def\thesection{A}
\def\theequation{A.\arabic{equation}}}

\newcommand{\mylabel}[1]{\def\nummer{#1}
\label{\nummer}}

\def\thesection{\arabic{section}}
\def\theequation{\arabic{section}.\arabic{equation}}

\advance\textheight-\headheight      
\advance\textheight-\headsep          
\advance\textheight by \topskip
\skip\footins 9pt plus 4pt minus 2pt
\begin{document}
\thispagestyle{empty}
\centerline{\bf M.~V.~Keldysh Institute of Applied Mathematics} 
\centerline{\bf Russian Academy of Sciences}
\vskip 2.5cm
\centerline{\bf R. Denk, R. Mennicken,  and L. R. Volevich}
\vskip 1.7cm
\centerline{\bf ON ELLIPTIC OPERATOR PENCILS}
\vskip 0.3cm
\centerline{\bf  WITH GENERAL BOUNDARY CONDITIONS} 
\vskip 2cm
\centerline{Preprint No. 37}
\vskip 5 cm
\centerline{\bf Moscow 1999}

\newpage
\thispagestyle{empty}
\centerline{Abstract}

R. Denk, R. Mennicken, and L. Volevich\footnote{Supported in part 
by the Deutsche Forschungsgemeinschaft and by Russian Foundation of 
Fundamental Research, Grant 97-01-00541}. On elliptic operator pencils
with general boundary conditions.

\medskip
In this paper operator pencils $A(x,D,\lambda)$ are investigated which 
depend polynomially on the parameter $\lambda$ and act on a manifold
with boundary. The operator $A$ is assumed to satisfy the condition of 
$N$-ellipticity with parameter which is an ellipticity condition
formulated with the use of the Newton polygon. We consider general
boundary operators $B_1(x,D),\ldots,B_m(x,D)$ and define
$N$-ellipticity for the boundary value problem $(A,B_1,\ldots,B_m)$
analogously to  the Shapiro--Lopatinskii condition.
It is shown that the boundary value problem is $N$-elliptic if and
only if an a priori estimate holds, where the norms in the estimate
are again defined in terms of the Newton polygon. These results are
closely connected with singular perturbation theory and lead to
uniform estimates for  problems of
Vishik--Lyusternik  type  containing a small parameter.

\bigskip



\newpage

\mysection{1}{Introduction}

In this paper we consider an operator pencil depending polynomially on 
the complex parameter $\lambda$ and being of the form
\begin{equation}\mylabel{1-1}
A(x,D,\lambda) = A_{2m}(x,D)+\lambda A_{2m-1}(x,D)
+ \cdots + \lambda^{2m-2\mu} A_{2\mu}(x,D)\,,
\end{equation}
where $m$ and $\mu$ are integer numbers with $m>\mu>0$ and $A_{2\mu},\ldots,
A_{2m}$ are partial differential operators of the form
\begin{equation}\mylabel{1-2}
 A_j(x,D) = \sum_{|\alpha|\le j} a_{\alpha j}(x) D^\alpha \quad (j=2\mu,\,
2\mu+1, \ldots, 2m)\,.
\end{equation}
We assume that the pencil \eqnref{1-1} acts on a compact manifold $M$
with boundary ${\partial M}$; the coefficients $a_{\alpha j}$ in
\eqnref{1-2} are complex-valued. The manifold, its boundary and the
coefficients are assumed to be infinitely smooth. In \eqnref{1-2} and
in the following, we write $x=(x_1,\ldots,x_n)$ and use
 the standard multi-index notation
\[
  D^\alpha = D_1^{\alpha_1}\cdots D_n^{\alpha_n}\,,\quad 
D_j = -i \frac{\partial}{\partial x_j}\,,\quad |\alpha| =
\alpha_1 + \ldots + \alpha_n\,.
\]
The Dirichlet boundary value problem corresponding to the pencil
\eqnref{1-1} was studied in detail in the paper \cite{dmv2}. The aim
of the present paper is to obtain an a priori estimate for general
boundary value problems connected with \eqnref{1-1}. So we assume that 
we have  boundary operators, for simplicity independent of the
complex parameter $\lambda$, of the form
\begin{equation}\mylabel{1-4}
 B_j(x,D) = \sum_{|\beta|\le m_j} b_{j\beta}(x) D^\beta \quad
 (j=1,\ldots, m)\,,
\end{equation}
where the numbering is chosen such that for the orders of the
operators $B_j$ we have $m_1\le m_2 \le \ldots \le m_m$. Additionally, 
we assume that 
\begin{equation}\mylabel{1-4a}
m_\mu < m_{\mu+1}\,.
\end{equation}
 The coefficients of $B_j$ are
supposed to be defined in $\overline M$ and to be infinitely smooth.

We will show that an a priori estimate holds if the boundary value
problem $(A,B_1,\ldots,B_m)$ satisfies the condition of
$N$-ellipticity with parameter which will be explained in Section
2. Moreover, we will prove that $N$-ellipticity is also necessary for
the validity of the a priori estimate (see Section 6).

The principal symbol $A^{(0)}(x,\xi,\lambda)$ of \eqnref{1-1} is
defined as
\begin{equation}\mylabel{1-5}
A^{(0)}(x,\xi,\lambda) := A_{2m}^{(0)}(x,\xi) + \lambda A_{2m-1}^{(0)}(x,
\xi) + \ldots + \lambda^{2m-2\mu}A_{2\mu}^{(0)}(x,\xi)\,,
\end{equation}
where
\begin{equation}\mylabel{1-6}
A_j^{(0)}(x,\xi) := \sum_{|\alpha|=j} a_{\alpha j}(x)\xi^\alpha\quad (j=2\mu,
\ldots,2m)
\end{equation}
stands for the principal symbol of $A_j$. In \eqnref{1-6} we have set 
$\xi^\alpha = \xi_1^{\alpha_1}\cdots
\xi_n^{\alpha_n}$ for $\xi=(\xi_1,\ldots,\xi_n)$.
 The principal symbols \eqnref{1-5} and \eqnref{1-6} are 
invariant under change of coordinates and thus globally defined on the 
cotangent bundle $T^*M\backslash\{0\}$. The principal symbols
$B_j^{(0)}$ of the boundary operators $B_j$ are defined analogously.

In \cite{dmv2} the Newton polygon approach was used to formulate and
prove an a priori estimate for the Dirichlet boundary value problem.
 This method (which was also applied to
Douglis--Nirenberg systems in \cite{dmv}) turns out to be suitable for 
general boundary conditions, too. The concept of the Newton polygon
makes it possible to define the general notion of $N$-ellipticity with
parameter which is a generalization of the
classical definition of ellipticity with parameter given by Agmon
\cite{agmon} and Agranovich--Vishik \cite{agranovich-vishik}. For the
connection to $N$-parabolic problems and Douglis--Nirenberg systems,
the reader is referred to \cite{dmv2}, Section 1.

Replacing in \eqnref{1-1} $\lambda$ by $\epsilon^{-1}$, we obtain a
problem of singular perturbation theory as it was studied, for
instance, by Vishik and Lyusternik \cite{vishik-lyusternik}. The a
priori estimate stated below in Section 4 corresponds to a uniform
(with respect to $\epsilon$) estimate in the Vishik--Lyusternik theory 
(see also \cite{frank}, \cite{nazarov}). We will show this close
connection in the Appendix.

\mysection{2}{The
  Shapiro--Lopatinskii condition}

As the manifold $M$ is compact we may fix a finite number of
coordinate systems. Locally in each of these coordinate systems the
operator pencil $A(x,D,\lambda)$ is of the form \eqnref{1-1} and acts
in $\R^n$. 
We can suppose without loss of generality that the coefficients 
of $A(x,D,\lambda)$ are (in local coordinates) of the form
\begin{equation}\mylabel{2-1}
a_{\alpha j}(x)= a_{\alpha j}+a'_{\alpha j}(x),\quad a'_{\alpha j}
\in {\cal D}(\R^n)  \,.
\end{equation}

\begin{definition}\mylabel{2.00}{\rm Let $x^0\in \overline M$ be fixed. The
 interior symbol
$A(x^0,\xi,\lambda)$ is called $N$-elliptic with parameter in
$[0,\infty)$ at $x^0$
(cf. \cite{dmv2}) if the estimate
\begin{equation}\mylabel{2-2}
| A^{(0)}(x^0,\xi,\lambda)| \ge C |\xi|^{2\mu}\,(\lambda + |\xi|)^{2m-2\mu}
\quad \big(\xi\in\R^n, \lambda\in [0,\infty)\big)
\end{equation}
holds with a constant $C$ which does not depend on $\xi$ or $\lambda$. If 
this is true for every $x^0\in\overline M$, the symbol
$A(x,\xi,\lambda)$ and the operator $A(x,D,\lambda)$ are called
$N$-elliptic with parameter in $[0,\infty)$.}
\end{definition}

By continuity and compactness, for an $N$-ellliptic operator the
constant $C$  in
\eqnref{2-2} can be chosen independently of $x^0$. 

Now we shall define the analogue of the Shapiro--Lopatinskii condition 
for our problem. For this,
we fix a point $x^0 \in {\partial M}$ and a coordinate system in the 
neighbourhood of $x^0$ such that in this system locally the
 boundary ${\partial M}$
is given by the equation $x_n=0$. 
We use in $\R^n_+ := \{ x\in\R^n: x_n>0\}
$ the coordinates $x=(x',x_n)$ and the dual coordinates
$\xi=(\xi',\xi_n)$. If $A$ is $N$-elliptic with parameter,
   it follows from \cite{dmv2}, Lemma 3.2, that for every $x^0\in
   \overline M$ we have
\begin{equation}\mylabel{2-5}
 A^{(0)}(x^0, \xi,\lambda)\not=0 \quad \big(\xi\in\R^n\backslash\{0\},\;
\lambda\in[0,\infty)\big)\,.
\end{equation}
In the case $n>2$ this implies that $A^{(0)}$, considered as a
polynomial in $\xi_n$, has exactly $m$ roots with positive imaginary
part for every $\xi'\not=0$.
 In the case $n=2$ this is an additional condition which we
assume to hold in the following. Similar considerations hold for 
$A_{2\mu}^{(0)}$.

Let $A$ be $N$-elliptic with parameter in $[0,\infty)$, fix
$x^0\in{\partial M}$ and write $A$ in local coordinates corresponding to
$x^0$ as considered above. Then we define 
 the polynomial in $\tau\in\C$ 
\begin{equation}\mylabel{2-3}
Q(x^0,\tau)=\tau^{-2\mu}A^{(0)}(x^0,0,\tau,1)   \,.  
\end{equation}

\begin{definition}\mylabel{2.0}{\rm The operator $A(x,D,\lambda)$
 degenerates regularly at the boundary ${\partial M}$ if for every
 $x^0\in{\partial M}$ the polynomial \eqnref{2-3}
 has exactly $m-\mu$ roots in the upper half-plane of the complex 
plane.}
\end{definition}

\begin{remark}\mylabel{2.01}{\rm 
a) It is easily seen that if for a fixed $x^0\in \partial M$
 and a fixed 
coordinate system polynomial \eqnref{2-3} has $m-\mu$
 roots in the upper half-plane, 
then this polynomial has this property for arbitrary $x^0\in \partial
 M$  and 
for an arbitrary coordinate system.
 
\noindent b) The condition of regular degeneration has its direct
counterpart in the theory of singular perturbations (see,
e.g., \cite{vishik-lyusternik}, Section 6). 

\noindent c) Some examples where the condition of regular degeneration 
(Definition \ref{2.0})
holds automatically can be found in \cite{dmv2}, Remark 3.4.
}
\end{remark}

If $A$ is $N$-elliptic with parameter in $[0,\infty)$, then for any
fixed $x^0\in\overline M$  and
$\xi'\in\R^{n-1}\backslash\{0\}$, we see from \eqnref{2-5} that we
can factorize the principal
symbol $A^{(0)}(x^0,\xi,\lambda)$ in the form
\[ A^{(0)}(x^0,\xi,\lambda) = A^{(0)}_+ (x^0,\xi,\lambda)\;
A^{(0)}_- (x^0,\xi,\lambda)\,.\]
Here  
\begin{equation}\mylabel{2-4}
   A^{(0)}_+ (x^0,\xi',\tau,\lambda) := \prod_{j=1}^m (\tau -
  \tau_{j}(x^0,\xi',\lambda)) \,,
  \end{equation}
  where $\tau_{1},\ldots,\tau_{m}$ are the zeros of 
$A^{(0)}$ with positive
  imaginary part.

Now let $x^0\in{\partial M}$ and 
denote by $B_j'(x^0,\xi',\xi_n,\lambda)$ the remainder of
$B_j^{(0)}(x^0,\xi)$ after division by $A^{(0)}_+(x^0,\xi,\lambda)$,
where all polynomials are considered as polynomials in $\xi_n$. We
write $B_j'$ in the form
\begin{equation}\mylabel{2-4a}
  B_j'(x^0,\xi',\xi_n,\lambda) = \sum_{k=1}^m b_{jk}(x^0,\xi',\lambda) 
  \xi_n^{k-1}\,.
\end{equation}
and define the Lopatinskii determinant by
\begin{equation}\mylabel{2-4b}
  \mbox{Lop}(x^0,\xi',\lambda) := \det\Big(
  b_{jk}(x^0,\xi',\lambda)\Big)_{j,k=1,\ldots,m}\,.
\end{equation}
Then the condition 
\begin{equation}\mylabel{2-4c}
\mbox{Lop}(x^0,\xi',\lambda) \not=0
\end{equation}
means that the polynomials $B_j^{(0)}(x^0,\xi',\cdot)$ are linearly
independent modulo $A^{(0)}_+(x^0,\xi',\cdot,\lambda)$.
It is well-known that condition \eqnref{2-4c} is satisfied if and only 
if the ordinary differential equation on the half-line
\begin{eqnarray}
A(\xi', D_t,\lambda)\, w(t)  & = &0 \quad\quad (t>0)\,,
\mylabel{eq2-1a}\\
B_k(\xi',D_t)\, w(t)|_{t=0} & = & h_k  \;\quad(k=1,\ldots,m)\,,
\mylabel{eq2-2a}\\
w(t) & \to & 0 \quad\quad (t \to + \infty)\,,\nonumber
\end{eqnarray}
is uniquely solvable for every $(h_1,\ldots,h_m)\in\C^m$.
Here $D_t$ stands for $-i\frac\partial{\partial t}$.

\begin{definition}\mylabel{2.1}{\rm Let $A$ satisfy the regular
    degeneration condition. Then the boundary problem
    $(A,B_1,\ldots, B_m)$ is called $N$-elliptic with parameter 
    $\lambda\in [0,\infty)$ if the following conditions hold:

\vskip 0.5em
\noindent a) The interior symbol $A(x,\xi,\lambda)$ is $N$-elliptic with 
  parameter in $[0,\infty)$ in the sense of Definition \ref{2.00}.

\vskip 0.5em
\noindent b) For every fixed $x^0\in{\partial M}$, every $\xi'\not=0$
 and every   $\lambda\in [0,\infty)$ the polynomials
  $(B_j^{(0)}(x^0,\xi',\cdot))_{j=1,\ldots,m}$ are linearly independent modulo
  $A^{(0)}_+ (x^0,\xi',\cdot,\lambda)$, i.e. \eqnref{2-4c} holds.

\vskip 0.5em
\noindent c) For every fixed $x^0 \in {\partial M}$, the
 boundary problem 
\[ (A_{2\mu}^{(0)}(x^0,D), B_1(x^0,D),\ldots,
  B_{\mu} (x^0,D))\]
 fulfills the Shapiro--Lopatinskii condition, 
i.e. 
  $(B_j^{(0)}(x^0,\xi))_{j=1,\ldots,\mu}$ are linearly independent
  modulo $(A_{2\mu}^{(0)})_+ (x^0,\xi)$. Here 
 $(A_{2\mu}^{(0)})_+ $ is defined in analogy
  to \eqnref{2-4} with $A$ replaced by $A_{2\mu}$.

\vskip 0.5em
\noindent d) Let $Q_+(x^0,\tau) := 
  \prod_{j=\mu+1}^m (\tau-\tau_j^1(x^0))$ where
  $\tau^1_{\mu+1}, \ldots, \tau^1_m$ denote the zeros of $Q(x^0,\tau)$ with
  positive imaginary part. Then the polynomials
  $(B_j^{(0)}(x^0,0,\tau))_{j=\mu+1,\ldots, m}$ are linearly independent 
  modulo $Q_+(x^0,\tau)$ for every $x^0\in{\partial M}$.

}
\end{definition}

\begin{remark}
{\rm a) Note that  the degree of  $B_j^{(0)}(x^0,0,\cdot)$
  is $m_j$ which may be greater than $2m-2\mu$.

\noindent b) Note that condition b) in Definition \ref{2.1} differs
from the Agmon--Agranovich--Vishik condition of ellipticity with
parameter. If the symbols $A(x,\xi,\lambda)$ and $B_j(x,\xi)$ are
homogeneous with respect to $(\xi,\lambda)$, the
Agmon--Agranovich--Vishik condition means that
\begin{equation}\mylabel{2-4d}
  \mbox{Lop}(x^0,\xi',\lambda) \ne 0 
  \mbox{ for } |\xi'|^2+\lambda^2=1,\; \lambda \ge 0\,.
\end{equation}
In particular, in this case the inequality \eqnref{2-4c} holds for
$\lambda=1$ and $\xi'=0$. In the case of $N$-ellipticity, however, the 
Lopatinskii determinant is in general not defined for $\xi'=0$ and may 
tend to zero as $\xi'\to0$.

\noindent c) Taking in \ref{2.1} b) $\lambda=0$ and $|\xi'|=1$, we
obtain the standard Shapiro--Lopatinskii condition for the boundary
value problem $(A_{2m},B_1,\ldots,B_m)$.

\noindent d) In the Vishik--Lyusternik theory, the analogue of
condition \ref{2.1} d) leads to the existence of solutions of boundary 
layer type (see \cite{vishik-lyusternik}, Section 6).
}
\end{remark}

\mysection{3}{The basic ODE estimate}

In a first step we consider the 
model problem in the half space.
Let $(A,B_1,\ldots,B_m)$ be of the form
\eqnref{1-1}, \eqnref{1-4} 
and acting in $\R^n_+$. We suppose
that $A$ is homogeneous in $(\xi,\lambda)$, i.e. has the form
\begin{equation}\mylabel{2-1a}
A(\xi,\lambda) = A_{2m}(\xi) + \lambda A_{2m-1}(\xi) +\ldots
+\lambda^{2m-2\mu} A_{2\mu}(\xi)\,,
\end{equation}
where $A_j(\xi)$ is a homogeneous polynomial in $\xi$ of degree $j$.
Similarly we assume that $B_j$ is given by
\begin{equation}\mylabel{2-1b}
B_j(\xi) = \sum_{|\beta|= m_j} b_{j\beta} \xi^\beta \quad
 (j=1,\ldots, m)\,.
\end{equation}

For fixed $\lambda\ge 0$ and $\xi'\in 
\R^{n-1}\backslash \{0\}$ we investigate the boundary problem 

\begin{eqnarray}
A(\xi', D_t,\lambda)\, w_j(t)  & = &0 \quad\quad (t>0)\,,
\mylabel{eq2-1}\\
B_k(\xi',D_t)\, w_j(t)|_{t=0} & = & \delta_{jk}  \quad(k=1,\ldots,m)\,,
\mylabel{eq2-2}\\
w_j(t) & \to & 0 \quad\quad (t \to + \infty)\,.\nonumber
\end{eqnarray}

In \cite{dmv2}, the following lemma on the roots of the polynomial
$A(\xi',\cdot,\lambda)$ is proved.

\begin{lemma}\mylabel{2.2a} Let the polynomial $A(\xi,\lambda)$ in
\eqnref{2-1a} be $N$-elliptic with parameter in $[0,\infty)$ and
assume that $A$ degenerates regularly.
 Then, with a suitable numbering of the roots $\tau_j(\xi',
\lambda)$ of $A(\xi',\tau,\lambda)$ with positive imaginary part, we 
have:\\
{\rm (i)} Let $S(\xi')= \{ \tau_1^0(\xi'),\ldots,\tau_\mu^0(\xi')\}$ be the 
set of all zeros of $A_{2\mu}(\xi',\tau)$ with positive imaginary
part. Then for all $r>0$ there exists a $\lambda_0>0$ such that the
distance between the sets $\{
\tau_1(\xi',\lambda),\ldots,\tau_\mu(\xi',\lambda)\}$ and $S(\xi')$ is 
less than $r$ for all $\xi'$ with $|\xi'|=1$ and all
$\lambda\ge\lambda_0$.\\
{\rm (ii)}  Let $\tau_{\mu+1}^1,\ldots,\tau_m^1$ be the roots of the polynomial
$Q(\tau)$ {\rm (}cf. \eqnref{2-3}{\rm)}
 with positive imaginary part. Then
 \begin{equation}\mylabel{eq2-16}
\tau_j(\xi',\lambda) = \lambda\tau_j^1 + \tilde\tau_j^1(\xi',\lambda)
\quad(j=\mu+1,\ldots,m)\,,
\end{equation}
and there exist  constants $K_j$ and $\lambda_1$, independent of
$\xi'$ and $\lambda$,
such that for $\lambda\ge \lambda_1$ the inequality
\begin{equation}\mylabel{eq2-17}
|\tilde\tau_j^1(\xi',\lambda)|\le K_j |\xi'|^{\frac 1{k_1}}\,
\lambda^{1-\frac 1{k_1}}\quad  (|\xi'|\le\lambda)
\end{equation}
holds, where $k_1$ is the maximal multiplicity of the roots of $Q(\tau)$.
\end{lemma}

\begin{theorem}\mylabel{2.2} Assume that the operator
$(A,B_1,\ldots,B_m)$ is of the form
\eqnref{2-1a}--\eqnref{2-1b}. Assume that condition \eqnref{1-4a}
holds and that $A$ degenerates regularly at the
boundary {\rm (}cf. Definition {\rm \ref{2.0}}{\rm)} 
 and $(A,B_1,\ldots,B_m)$ is $N$-elliptic with parameter in $\R^n_+$ in the
  sense of Definition {\rm \ref{2.1}}. Then
for every $\xi'\in\R^{n-1}\backslash\{0\}$
 and $\lambda\in [0,\infty)$ the ordinary differential equation 
\eqnref{eq2-1}--\eqnref{eq2-2} has a unique solution $w_j(t,\xi',\lambda)$,
and the estimate
\begin{eqnarray}
& & \hspace*{-5em}
\| D_t^l w_j(\cdot,\xi',\lambda)\|_{L_2(\R_+)}\le  \nonumber\\
& & \hspace*{-4em}\le C \left\{
\def\arraystretch{1.2}
\begin{array}{lll}
|\xi'|^{l-m_j -\frac 1 2}, & j\le \mu,&l\le m_{\mu+1},\\
|\xi'|^{m_{\mu+1}-m_j}(\lambda + |\xi'|)^{l-m_{\mu+1}- \frac 1 2},
 & j \le \mu,& l> m_{\mu+1},\\
|\xi'|^{l-m_\mu-\frac 1 2}(\lambda + |\xi'|)^{m_\mu-m_j},
 & j > \mu,& l\le m_\mu,\\
(\lambda + |\xi'|)^{l-m_j-\frac 1 2}, & j > \mu,& l> m_\mu,
\end{array}\right.
\def\arraystretch{1}\mylabel{eq2-3}
\end{eqnarray}
holds with a constant $C$ not depending on $\xi'$ and $\lambda$.
\end{theorem}

\begin{proof}
The existence and the uniqueness of the solution follows immediately
from conditions a) and b) in Definition \ref{2.1}. \mbox{}From the
homogeneity of the symbols  and from the uniqueness of the solution we see 
that
\begin{equation}\mylabel{eq2-4}
  w_j(t,\xi',\lambda) = r^{-m_j} w_j\Big(rt,\frac{\xi'}{r}, 
  \frac\lambda r\Big)
\end{equation}
holds for every $r>0$. If we set 
$r=|\xi'|$ and $\omega' = \frac{\xi'}{|\xi'|}$ we obtain
\begin{equation}\mylabel{eq2-5}
\| D_t^lw_j(\cdot,\xi',\lambda)\|_{L_2(\R_+)} = |\xi'|^{
  l-m_j-\frac 1 2}
\Big\| D_t^l w_j\Big(\cdot,\omega',  \frac\lambda {|\xi'|}\Big)
\Big\|_{L_2(\R_+)}\,.
\end{equation}
The theorem will be proved if we show that for $|\omega'|=1$ we have
\def\arraystretch{1.2}
\begin{equation}\mylabel{eq2-6}
\|(D_t^l w_j)(\cdot,\omega',\Lambda)\|_{L_2(\R_+)} \le
 \left\{\begin{array}{ll}
 C\,, & j\le \mu\,,\; l\le m_{\mu+1}\,,\\
 C\,\Lambda^{l-m_{\mu+1}- \frac 1 2}, & j\le \mu\,,\; l> m_{\mu+1}\,, \\
 C\,\Lambda^{m_\mu-m_j}, & j> \mu\,,\; l\le m_\mu\,,\\
 C\,\Lambda^{l-m_j-\frac 1 2}, & j> \mu\,,\; l> m_\mu\,,\\
 \end{array}\right.
\end{equation}
\def\arraystretch{1}
for $\Lambda\ge 1$ and that the left-hand side is bounded by a
constant for $\Lambda\le 1$. 

The boundedness for $\Lambda\le 1$ follows easily from conditions a)
and b) of Definition \ref{2.1}. We have to consider the case of large
$\Lambda$. 

To find an estimate in this case, we represent the solution in a form
suggested in a paper of Frank \cite{frank}. This representation is
different from the formula used for the Dirichlet boundary value
problem in \cite{dmv2}  and allows us to separate
the two parts of the zeros of the polynomial $A_+(\omega',\tau,\Lambda)$ 
in a more adequate form.

Due to Lemma \ref{2.2a}, the roots of this polynomial consist of two 
groups, the first group, denoted by $\{\tau_1(\omega',\Lambda),\ldots,
\tau_\mu(\omega',\Lambda)\}$, being bounded for $\Lambda\to\infty$, the 
other group, denoted by $\{\tau_{\mu+1}(\omega',\Lambda),\ldots,
\tau_m(\omega',\Lambda)\}$, being of order $\Lambda$ for $\Lambda \to
\infty$. 

We define
\begin{equation}\mylabel{eq2-7}
A_1(\omega',\tau,\Lambda) := \prod_{j=1}^\mu (\tau
-\tau_j(\omega',\Lambda))\,. 
\end{equation}
Let $\gamma^{(1)}$ be a contour in the upper half of the complex
plane enclosing the zeros $\tau_1,\ldots, \tau_\mu$. \mbox{}From Lemma
\ref{2.2a}  we 
see that $\gamma^{(1)}$ can be chosen independently of $\omega'$ and
$\Lambda$ for all $|\omega'|=1$ and $\Lambda\ge \Lambda_0$. 

\mbox{}From the same lemma we see that $A_1(\omega',\tau,\Lambda)\to
(A_{2\mu})_+ (\omega',\tau)$ as $\Lambda\to\infty$. Therefore we
obtain from condition c) in Definition \ref{2.1} that for all
$|\omega'|=1$ and $\Lambda\ge \Lambda_0$ the polynomials $\{
B_j(\omega', \tau)\}_{j=1,\ldots,\mu}$ are independent modulo
$A_1(\omega', \tau,\Lambda)$. Thus there exist polynomials (with
respect to $\tau$) $N_j(\omega',\tau,\Lambda)$, depending continuously 
on $(\omega',\Lambda)$, such that
\begin{equation}\mylabel{eq2-8}
\frac{1}{2\pi i} \int_{\gamma^{(1)}} \frac{B_k(\omega',\tau)
  N_j(\omega', \tau,\Lambda)}{A_1(\omega',\tau,\Lambda)} d\tau =
\delta_{kj} \quad (k,j=1,\ldots,\mu)\,.
\end{equation}
(For the construction of $N_j$ cf.,
e.g., \cite{agmon-douglis-nirenberg}, p. 634.)

Analogously, we define
\begin{equation}\mylabel{eq2-9}
A_2(\omega',\tau,\Lambda) := \prod_{j=\mu+1}^m (\tau
-\tau_j(\omega',\Lambda))\,. 
\end{equation}
Let $\tilde\gamma^{(2)}(\omega',\Lambda)$ be a contour in the 
upper half of the
complex plane enclosing the zeros $\tau_{\mu+1}(\omega',\Lambda),
\ldots, \tau_{m}(\omega',\Lambda)$. \mbox{}From Lemma \ref{2.2a} we know that this
contour is of order $\Lambda$ for $\Lambda\to\infty$. Therefore we may 
fix a contour $\gamma^{(2)}$, independent of $\omega'$ and $\Lambda$
such that $\gamma^{(2)}$ encloses all values $\tau_j  /\Lambda$
with $j=\mu+1,\ldots,m$. We also remark that due to the regular
degeneration  we may choose $\gamma^{(2)}$ with a positive
distance to the real axis (cf. also \eqnref{eq2-16}).

\mbox{}From condition d) in \ref{2.1} we know that $\{
B_j(0,\tau)\}_{j=\mu+1,\ldots, m}$ is linearly independent modulo
  $Q_+(\tau)$. \mbox{}From Lemma \ref{2.2a} b) we know that 
\[ A_2\Big(\frac{\omega'}{\Lambda},\tau,1\Big) \to Q_+(\tau) \quad
(\Lambda \to \infty)\,.\]
Due to continuity, the polynomials $\{
B_j(\frac{\omega'}{\Lambda},\tau, 1)\}_{j=\mu+1,\ldots,m}$ are for
sufficiently large $\Lambda$ linearly independent modulo
$A_2(\frac{\omega'}{\Lambda}, \tau, 1)$. Therefore there exist
polynomials (in $\tau$) $N_j(\omega', \tau,\Lambda)$ for $j=\mu+1,
\ldots, m$, depending continuously on $\omega'$ and $\Lambda$, such
that
\begin{equation}\mylabel{eq2-10}
\frac{1}{2\pi i} \int_{\gamma^{(2)}} \frac{B_k(\frac{\omega'} \Lambda, 
  \tau) N_j(\omega', \tau,\Lambda)}{A_2(\frac{\omega'}\Lambda, \tau, 1 
  )} d\tau =
\delta_{kj} \quad (k,j=\mu+1,\ldots,m)\,.
\end{equation}

Now we need a lemma which will be proved below.

\begin{lemma}\mylabel{2.3}
 The solution $w_j(t,\omega',\Lambda)$ of the problem 
\eqnref{eq2-1}--\eqnref{eq2-2}
can be represented in the form
\begin{eqnarray}
w_j(t,\omega',\Lambda) & = &
\frac{1}{2\pi i} \int_{\gamma^{(1)}} \frac{M^{(1)}_j
(\omega', \tau,\Lambda)}{A_1(\omega',\tau,\Lambda)}
e^{it\tau}d\tau \nonumber \\
& + &\frac{1}{2\pi i} \int_{\gamma^{(2)}} \frac{M^{(2)}_j
(\omega', \tau,\Lambda)}{A_2(\frac{\omega'}{\Lambda},\tau,1)}
e^{it\Lambda\tau}d\tau \mylabel{2-6}
\end{eqnarray}
where for $|\tau|=O(1)$ and $|\omega'|=1$ we have
\[ M_j^{(1)}(\omega',\tau,\Lambda) \le \left\{ 
   \begin{array}{ll}
    C\,, &  j\le \mu\,,\\
    C\,  \Lambda^{m_{\mu}-{m_j}}\,, &  j> \mu\,,
    \end{array}\right. \]
and
\[ M_j^{(2)}(\omega',\tau,\Lambda) \le\left\{ 
    \begin{array}{ll}
    C\,\Lambda^{-m_{\mu+1}}\,,& j \le \mu\,,\\
    C\,\Lambda^{-m_j}\,,         & j> \mu\,,
    \end{array}\right. \]
\end{lemma}

As a direct corallary of the lemma we obtain
\[ \|(D_t^l w_j)(\cdot,\omega',\Lambda)\|_{L_2(\R_+)} \le
\left\{ \begin{array}{ll}
O(1)+O(\Lambda^{l-m_{\mu+1}-\frac{1}{2}})\,,& j\le \mu\,,\\
O(\Lambda^{m_{\mu}-m_j})+O(\Lambda^{l-m_j-\frac {1}{2}})\,,
& j > \mu\,.
\end{array}\right. \]
The estimate \eqnref{eq2-6} trivially follows from these relations.
\end{proof}

\begin{proofof}{Proof of Lemma {\rm \ref{2.3}}}
Let $w(t,\omega',\Lambda)$ be a solution of the problem
\eqnref{eq2-1}--\eqnref{eq2-2}
 with $\delta_{jk}$ replaced by $\phi=(\phi_1,\dots,\phi_m)\in\C^m$.
We seek the solution in the form
\begin{eqnarray}
w(t,\omega',\Lambda)&=&\sum_{k=1}^{\mu}\psi_k(\omega',\Lambda)
\frac{1}{2\pi i} \int_{\gamma^{(1)}} \frac{N_k
(\omega', \tau,\Lambda)}{A_1(\omega',\tau,\Lambda)}
e^{it\tau}d\tau\nonumber\\
& & \hspace*{-3em}+  \sum_{k=\mu+1}^m \psi_k(\omega',\Lambda)
\frac{1}{2\pi i} \int_{\gamma^{(2)}} \frac{N_k
(\omega', \tau,\Lambda)}{A_2(\frac{\omega'}{\Lambda},\tau,1)}
e^{it\Lambda\tau}d\tau\mylabel{2-7}
\end{eqnarray}
where the functions $\psi_k$ still have to be found.

Applying the  boundary  operator  $B_l(\xi',D_t)$  to both
sides of \eqnref{2-7} and taking $t=0$ we obtain the following system
for the unknown functions $\psi_k(\omega',\Lambda)$:

\begin{eqnarray}
\psi_l(\omega',\Lambda) +\Lambda^{m_l} \sum_{k=\mu+1}^m \psi_k(\omega',\Lambda)
h_{lk} (\omega',\Lambda) & = & \phi_l \nonumber\\
& & \hspace*{-4em} (l=1,\ldots,\mu)\,,\mylabel{2-12}\\
\sum_{k=1}^\mu \psi_k(\omega',\Lambda)h_{lk} (\omega',\Lambda)
 + \Lambda^{m_l}
\psi_l(\omega',\Lambda) & = &  \phi_l\nonumber\\
& & \hspace*{-4em} (l=\mu+1,\ldots,m)\,.\mylabel{2-13}
\end{eqnarray}
Here we have set

\begin{eqnarray}
\lefteqn{h_{lk} (\omega',\Lambda)  =  \frac{1}{2\pi i} \int_{\gamma^{(2)}}
\frac{ B_l(\frac{\omega'}{\Lambda},
\tau) N_k(\omega',\Lambda,\tau)}{ A_2(\frac
  {\omega'}{\Lambda}, \tau, 1)}\; d \tau }\nonumber \\
& & \hspace*{7em} (l=1,\ldots,\mu; \;
k=\mu+1,\ldots, m)\,,\mylabel{2-14}\\[1em]
\lefteqn{ h_{lk} (\omega',\Lambda)  =  \frac{1}{2\pi i} \int_{\gamma^{(1)}}
\frac{ B_l(\omega', \tau) N_k(\omega',\Lambda,\tau)}{ A_1(
  {\omega'}, \tau,\Lambda)}\; d \tau}\nonumber\\
& &  \hspace*{7em}\quad (l=\mu+1,\ldots,m; \;
k=1,\ldots, \mu)\,.\mylabel{2-15}
\end{eqnarray}
We remark that  we have used $B_l(\omega', \Lambda\tau)
= \Lambda^{m_l} B_l(\frac{\omega'}{\Lambda},\tau)$.

Now we  write $\psi=(\psi',\psi'')$,  where $\psi'$ consists of the
first $\mu$ components of the vector $\psi$, and $\psi''$ consists
of the other  $m-\mu$ components.  In the same way we write $\phi=
(\phi',\phi'')$. In these notations the system  
\eqnref{2-12}--\eqnref{2-13}  can    
be rewritten in the form   
\def\arraystretch{1.2}
\[ \begin{array}{rrrrr}
\psi' &+&\Delta_1 H_{12}\psi'' &= &\phi'\,,\\
H_{21}\psi'&+&\Delta_2\psi''   &= &\phi''\,,  
\end{array}\] 
\def\arraystretch{1}
where we use the notation 
\[ \Delta_1 := \left( \begin{array}{ccc} \Lambda^{m_1} & & \\ &
    \ddots & \\ & & \Lambda^{m_\mu}\end{array}\right),\quad
  \Delta_2 :=  \left( \begin{array}{ccc} \Lambda^{m_{\mu+1}} & & \\ &
    \ddots & \\ & & \Lambda^{m_m}\end{array}\right) \]
and  
\[ H_{12} := \Big( h_{lk}\Big)_{\textstyle{l=1,\ldots,\mu \atop
  k=\mu+1,\ldots,m}}\,,\quad
H_{21} := \Big( h_{lk}\Big)_{\textstyle{l=\mu+1,\ldots,m \atop
  k=1,\ldots,\mu}}\,.\]

If me multiply the second equation  by  the  matrix  $\Delta_1
H_{12}\Delta^{-1}_2$ from the left
 and  subtract it from  the  first equation we
obtain
\[ (I-\Delta_1 H_{12}\Delta^{-1}_2 H_{12})\psi'=\phi'-\Delta_1
H_{12}\Delta_2^{-1}\phi''\,.\]
In a similar way we obtain
\[ (I-\Delta^{-1}_2 H_{21}\Delta_1 H_{12})\psi''=-\Delta_2^{-1}
H_{21}\phi'+\Delta_2^{-1}\phi''\,.\]
The matrices in brackets  in  the  left-hand  sides  of  above
relations differ  from the  identity by matrices whose elements can be 
estimated by a constant times $\Lambda^{m_\mu-m_{\mu+1}}$. According
to \eqnref{1-4a}, their norms tend to zero as $\Lambda\to\infty$.
\mbox{}From
this it follows  that the  matrices in brackets  for  large $\Lambda$ have
inverses which we denote by $G_1$ and $G_2$, respectively. Then
we obtain
\def\arraystretch{1.2}
\[ \begin{array}{rrrrr} 
\psi' &=& G_1\phi' &-& G_1\Delta_1 H_{12}\Delta_2^{-1}\phi''\,,\\
\psi''&=&-G_2\Delta^{-1}_2 H_{21}\phi'&+&G_2\Delta^{-1}_2\phi''\,.
\end{array} \] 
\def\arraystretch{1} 
If we take $\phi=e_j\quad  (1\le j\le \mu)$, where $e_j$ stands for
the $j$-th unit vector,  and  denote by $e'_j$  the 
first $\mu$ components of $e_j$, we obtain
\[ \psi'_{(j)}=G_1e'_j,\quad \psi''_{(j)}=-G_2\Delta^{-1}_2 H_{21} 
e'_j\,.\]
In the same way if  $j>\mu$  and  $e''_j$  denotes the  components
$\mu+1,\dots,m$ of $e_j$, we obtain
\[ \psi'_{(j)}=-G_1\Delta_1 H_{12}\Lambda^{-m_j}e''_j,\quad
\psi_{(j)}''=G_2\Lambda^{-m_j}e''_j\,.\]
The statement  of  the  lemma  directly  follows  from   these
relations. 
\end{proofof}

\mysection{4}{A priori estimates}

Theorem \ref{2.2} is the key result for proving a priori
estimates. The norms used in these estimates are based on the Newton
polygon $N_{r,s}$ (cf. Fig. 1) defined for $r>s\ge 0$
 as the convex hull of the set
\[ \{ (0,0)\,,\; (0,r-s)\,,\;(s,r-s)\,,\;(r,0)\,\}\]

\begin{figure}[ht]
\setlength{\unitlength}{1mm}
\begin{center}
\begin{picture}(90,60)(-5,5)
\put(5,15){\vector(1,0){80}}
\put(5,15){\vector(0,1){45}}
\put(82,10){$i$}
\put(0,55){$k$}
\put(5,41){\line(1,0){35}}
\put(40,41){\line(1,-1){25.8}}
\put(4.4,40.6){{$\scriptscriptstyle \bullet$}}
\put(39.4,40.6){{$\scriptscriptstyle \bullet$}}
\put(65,14.6){{$\scriptscriptstyle \bullet$}}
\put(40,14){\line(0,1){2}}
\put(65,10){$r$}
\put(39,10){$s$}
\put(-6,40){$r-s$}
\put(5,2){\parbox{60mm}{\begin{center}
{\small Fig. 1. The Newton polygon $N_{r,s}$.}\end{center}}}
\end{picture}
\end{center}
\end{figure}

The  weight function $\Xi_{r,s}(\xi,\lambda)$ is defined by
\begin{equation}\mylabel{3-1}
\Xi_{r,s}(\xi,\lambda) := \sum_{(i,k)\in N_{r,s}} |\xi|^i\,|\lambda|^k\,,
\end{equation}
where the summation on the right-hand side is extended over all 
integer points of $N_{r,s}$. For a discussion of general Newton
polygons we refer the reader to \cite{dmv}, \cite{dmv2},
\cite{gindikin-volevich}.

It is easily seen that we have the equivalence
\begin{equation}\mylabel{3-2}
 \Xi_{r,s} (\xi,\lambda) \approx (1+|\xi|)^{s} (\lambda + 
|\xi|)^{r-s}\,.
\end{equation}
The sign $\approx$ means that the quotient of the left-hand and the
right-hand side is bounded from below and from above by positive
constants independent of $\xi$ and $\lambda$. Taking the right-hand
side of \eqnref{3-2} as a definition, we may define $\Xi_{r,s}$ for
every $r,s\in\R$. The Sobolev space $H^{(r,s)}(\R^n) :=
H^{\Xi_{r,s}}(\R^n)$ is defined
as
\[ \{ u\in {\cal S}'(\R^n): \Xi_{r,s}(\xi,\lambda) Fu(\xi)\in
L_2(\R^n) \} \]
with the norm
\begin{equation}\mylabel{3-2a}
 \|u\|_{(r,s),\R^n} := \| u\|_{\Xi_{r,s},\R^n} := \| F^{-1} \Xi_{r,s}
(\xi,\lambda) Fu(\xi)\|_{L_2(\R^n)}\,.
\end{equation}
Here $Fu$ stands for the Fourier transform of $u$ and ${\cal
  S}'(\R^n)$ denotes the space of all tempered distributions. The space
$H^{\Xi_{r,s}}(\R^{n-1})$ is defined analogously with the weight function
$\Xi_{r,s}(\xi',\lambda) := \Xi_{r,s}(\xi',0,\lambda)$. These spaces
can be defined on the half-space $\R^n_+$ in accordance with the
general theory of Sobolev spaces with weight functions as it can be
found, e.g.,  in \cite{volevich-paneah}. On the manifold $M$ and the
boundary ${\partial M}$, the spaces $H^{\Xi_{r,s}}(M)$ and
$H^{\Xi_{r,s}}({\partial M})$, respectively, are defined in the usual way,
using a partition of unity. 

In \cite{dmv2}, Section 2, Sobolev spaces connected with Newton
polygons were investigated in detail. In particular, for $r\in \N$ 
it was shown that 
$(\frac{\partial}{\partial\nu})^j: u\mapsto
(\frac{\partial}{\partial\nu})^j u|_{\partial M}$ acts continuously from  
$H^{\Xi_{r,s}}(M)$ to $H^{\Xi_{r,s}^{(-j-1/2)}}({\partial M})$ for 
$j=0,\ldots, r-1$. Here $\Xi_{r,s}^{(-j-1/2)}$ 
denotes the weight function corresponding to the Newton
polygon  which is constructed from 
$N_{r,s}$ by a shift of length $j+1/2$ to the left parallel to the
abscissa. 

While in \cite{dmv2} the basic Sobolev space was 
$H^{\Xi_{m,\mu}}$, we here consider more general a priori
estimates. In the following, we fix integer numbers 
\begin{equation}\mylabel{3-2b}
 r\ge m_m+1 \; \mbox{ and }\; 
m_\mu+1 \le s \le m_{\mu+1}
\end{equation}
 and consider the Newton polygon $N_{r,s}$, its weight function
$\Xi:=\Xi_{r,s}$ and the corresponding 
 Sobolev space. We remark that for the Dirichlet problem the values
 $r=m$ and $s=\mu$ used in \cite{dmv2} are included as an example.

Analogously to \eqnref{3-1}, we define the function
$\Phi=\Phi_{r,s}$  by
\begin{equation}\mylabel{3-3} 
\Phi(\xi,\lambda) := \sum_{i,k} |\xi|^i \lambda^k\,,
\end{equation}
where the sum runs over all integer points $(i,k)$ belonging to the
side of $N_{r,s}$ which is not parallel to one of the coordinate
lines. This means that we have
\begin{equation}\mylabel{3-4}
\Phi(\xi,\lambda) \approx |\xi|^{s} (\lambda + 
|\xi|)^{r-s}\,.
\end{equation}
By $\Phi^{(-l)}$ we again denote the corresponding function for the
shifted Newton polygon.
\mbox{}From Theorem \ref{2.2} we obtain the following estimate for the
fundamental solution $w_j$ defined in \eqnref{eq2-1}--\eqnref{eq2-2}:

\begin{lemma}\mylabel{3.1} For the solution $w_j(t,\xi',\lambda)$
  considered in Theorem {\rm \ref{2.2}} we have the estimate
\begin{equation}\mylabel{3-5}
\| D_t^l w_j(\cdot, \xi', \lambda)\|_{L_2(\R_+)} \le C \;
\frac{ \Phi^{(-m_j-1/2)}(\xi',\lambda)}{\Phi^{(-l)}(\xi',\lambda)}\,.
\end{equation}
\end{lemma}

\begin{proof}
To see this, we only have to remark that the right-hand side of
\eqnref{3-5} is equivalent to
\[  \left\{
\def\arraystretch{1.2}
\begin{array}{lll}
|\xi'|^{l-m_j -\frac 1 2}, & j\le \mu,&l\le s,\\
|\xi'|^{s-m_j- \frac 1 2}(\lambda + |\xi'|)^{l-s},
 & j \le \mu,& l> s,\\
|\xi'|^{l-s}(\lambda + |\xi'|)^{s-m_j-\frac 1 2},
 & j > \mu,& l\le s,\\
(\lambda + |\xi'|)^{l-m_j-\frac 1 2}, & j > \mu,& l> s.
\end{array}\right.
\def\arraystretch{1} \]
The first and fourth lines above coincide with the corresponding lines in the
right-hand side of \eqnref{eq2-3}. 
The ratio of the second line in  \eqnref{eq2-3} and the second
line above is equal to
\[ \Big(\frac {|\xi'|}{\lambda+|\xi'|}\Big)^{m_{\mu+1}-s+1/2}\,. \]
Respectively, the ratio of the third line in \eqnref{eq2-3}
 and the third line above
 is equal to 
\[ \Big(\frac {|\xi'|}{\lambda+|\xi|}\Big)^{s-m_{\mu}-1/2}\,.\]
Now our statement follows from \eqnref{3-2b}.
\end{proof}

\begin{theorem}\mylabel{3.2} Let $A(x,D,\lambda)$ be an operator
  pencil of the form \eqnref{1-1}, acting on the manifold $M$ with
  boundary ${\partial M}$. Let $B_j(x,D)$, $j=1,\ldots, m$, be boundary
  operators of the form \eqnref{1-4}. Assume that $A$ degenerates
  regularly at the boundary and that $(A,B_1,\ldots,
  B_m)$ 
is  $N$-elliptic with parameter in the sense of Definition
 {\rm \ref{2.1}}. Set $\Xi = \Xi_{r,s}$ with  $r$ and $s$ 
 satisfying \eqnref{3-2b}. For simplicity, assume that $r$ and $s$ are
  integers. 
   Then for $\lambda\ge\lambda_0$ there exists a
  constant $C=C(\lambda_0)$, independent of $u$ and $\lambda$, such
  that
\begin{eqnarray}
\| u \|_{\Xi,M} & \le & C\Big( 
\|A(x,D,\lambda)u\|_{(r-2m, s-2\mu),M} \nonumber\\
&  & \hspace*{-5em} + \;
\sum_{j=1}^m \| B_j(x,D) u\|_{\Xi^{(-m_j-1/2)},{\partial M}} 
 + \lambda^{r-s} \|u\|_{L_2(M)}\Big)\,.\mylabel{3-6}
\end{eqnarray}
\end{theorem}

\begin{proof} The proof of this theorem is similar to the proof of
  Theorem~5.6 in \cite{dmv2}; therefore we only indicate the main
  steps.

By the localization method (``freezing the coefficients''), it is
possible to reduce the proof to the proof of the corresponding results 
for model problems in $\R^n$ and $\R^n_+$. The case of the whole space 
$\R^n$ is quite elementary and needs only slight changes in comparison 
with \cite{dmv2}, Proposition 5.2. The key result is the a priori
estimate in the half space $\R^n_+$.

So we assume that $u\in H^\Xi(\R^n_+)$ is a solution of
\begin{eqnarray}
\mylabel{3-7} A(D,\lambda) u(x) & = & f \quad \mbox{ in }\R^n_+\,,\\
\mylabel{3-8} B_j(D)u |_{x_n=0} & = & g_j\quad (j=1,\ldots,m)
\mbox{ on }\R^{n-1}\,.
\end{eqnarray}
 Let $\psi \in C^{\infty}(\R^n)$ be a cut-off function, i.e. $\psi
(\xi)=1$ for $|\xi| \le 1$ and $\psi(\xi)=0$ for $|\xi| \ge 2$. 
As in \cite{dmv2} we write $u$ in the form
\begin{equation}\mylabel{3-9}   
  u = u_1+u_2+v=R\psi(D)Eu+R(1-\psi(D)) A^{-1}(D,\lambda)Ef+v \,.
\end{equation} 
Here we fixed an extension operator $E$ from $\R^n_+$ to $\R^n$, being 
continuous from $L_2(\R^n_+)$ to $L_2(\R^n)$ and from
$H^{\Xi}(\R^n_+)$ to $H^{\Xi}(\R^n)$; we define the distribution $Ef$
as $A(D,\lambda)Eu$.  By $R$ we denote the operator
of restriction of functions on $\R^n$ onto $\R^n_+$. In \eqnref{3-9},
the pseudodifferential operator (ps.d.o.) $\psi(D)$ in $\R^n$ is
defined as 
\[ \psi(D) := F^{-1} \psi(\xi) F\,.\]

It is easily seen that we have
\begin{equation}\mylabel{3-10}
\|u_1\|_{\Xi,\R^n_+} + \|u_2\|_{\Xi,\R^n_+}\le C\Big(
\|f\|_{(r-2m, s-2\mu),\R^n_+} 
+ \lambda^{r-s} \|u\|_{L_2(\R^n_+)}\Big)\,.
\end{equation}
We still have to estimate $v$ defined in \eqnref{3-9}. By definition,
$v$ is a solution of
\begin{eqnarray}
\mylabel{3-11} A(D,\lambda)\,v & = & 0 \quad \mbox{ in }\R^n_+\,,\\
\mylabel{3-12} B_j(D)v|_{x_n=0} & = & h_j\quad (j=1,\ldots,m)
\mbox{ on }\R^{n-1}\,,
\end{eqnarray}
where we set $h_j(x') :=  B_j(D) u(x',0)-B_j(D) u_1(x',0)-
B_j(D)u_2(x',0)$.

Now we use the fact that 
\[ \|v\|_{\Xi,\R^n_+} \approx 
\|v\|_{\Phi,\R^n_+}+ \lambda^{r-s} \|v\|_{L_2(\R^n_+)}\,.\]
As
\[  \lambda^{r-s} \|v\|_{L_2(\R^n_+)}\le 
\lambda^{r-s}\|u\|_{L_2(\R^n_+)}+ \|u_1\|_{\Xi,\R^n_+}+
\|u_2\|_{\Xi,\R^n_+}\,,\]
it is sufficient to estimate $\|v\|_{\Phi,\R^n_+}$.

It can be seen, using the binomial formula, that
\[\big( \Phi(\xi,\lambda)\big)^2 \approx \sum_{l=0}^r \big(
\xi_n^l\; \Xi^{(-l)}(\xi',\lambda)\big)^2 \,,\]
and therefore  the
semi-norm $\|v\|_{\Phi,\R^n_+}$ is equivalent to
\begin{equation}\mylabel{3-13}
\left[ \sum_{l=0}^{r} \int_0^\infty \| (D_n^l
  v)(\cdot,x_n)\|^2_{\Phi^{(-l)}, \R^{n-1}}dx_n\right]^{1/2}\,.
\end{equation}
Taking the partial Fourier transform $F'$ with respect to
$x'\in\R^{n-1}$, we obtain from \eqnref{3-11}--\eqnref{3-12} that
for $\xi'\not=0$ the function
$w:= F'v$ is a solution of
\begin{eqnarray}
\mylabel{3-14} A(\xi',D_n,\lambda) w(x_n) & = & 0 \,,\\
\mylabel{3-15} B_j(\xi',D_n)w(x_n)|_{x_n=0} & = & (F'h_j)(\xi')\,.
\end{eqnarray}
Due to Theorem \ref{2.2}, this solution is unique and given by
\begin{equation}\mylabel{3-16}
w = w(x_n,\xi',\lambda) = \sum_{j=1}^m w_j(x_n,\xi',\lambda) (F'h_j)(\xi')
\end{equation} 
with $w_j(x_n,\xi',\lambda)$ being the solution of
\eqnref{eq2-1}--\eqnref{eq2-2}. Now we can apply Lemma \ref{3.1} to
obtain
\begin{eqnarray*}
\lefteqn{\hspace*{-1.5cm} (\Phi^{(-l)}(\xi',\lambda))^2\int_0^{\infty}
|D_n^l w(x_n,\xi',\lambda)|^2\,dx_n }\\
&\le &
C\; \sum\Big|\Xi^{(-m_j-\frac 1 2)}(\xi',\lambda)(F'h_j)(\xi')\Big|^2 \,.
\end{eqnarray*}
Integrating this inequality with respect to $\xi'$ and using the norm
\eqnref{3-13}  we get 
the desired estimate for  $\|v\|_{\Phi,\R^n_+}$, which finishes the
proof of the theorem.
\end{proof}

\mysection{5}{The parametrix construction}

In this section, we will construct a right (rough) parametrix for the
operator $(A,B) = (A,B_1,\ldots,B_m)$, generalizing the result of 
\cite{dmv2}. We restrict ourselves to the construction of local
parametrices in $\R^n$ and $\R^n_+$; after this the definition of the
parametrix  on the manifold 
 is standard
(cf. also \cite{dmv2} for the Dirichlet problem). 

For the remainder of this section, we fix integer numbers $r$ and $s$ 
satisfying \eqnref{3-2b} and set
$\Xi=\Xi_{r,s}$. The following result is a slight generaliztion of
\cite{dmv2}, Proposition 6.1, which can be proved literally in the
same way.

\begin{lemma}\mylabel{5.1} Let $A(x,D,\lambda)$ in \eqnref{1-1}
 be $N$-elliptic in
  $\R^n$ with coefficients of the form \eqnref{2-1}. Then there exists 
  a bounded operator
\[ P_0: H^{(r-2m,s-2\mu)}(\R^n) \to H^{(r,s)}(\R^n)\]
such that 
\[ AP_0 = I + T \]
where $I$ denotes the identity operator in $H^{(r-2m,s-2\mu)}(\R^n)$
and 
\[ T : H^\Theta(\R^n) \to H^{(r-2m,s-2\mu)}(\R^n)\]
is bounded. Here we have set
\begin{equation}\mylabel{5-1}
\Theta(\xi,\lambda) := \Xi_{r-2m+1,s-2\mu+1}(\xi,\lambda)\; \Big[ =
(1+|\xi|)\; \Xi_{r-2m,s-2\mu}(\xi,\lambda)\Big]\,.
\end{equation}
\end{lemma}

Throughout this section, by a bounded operator we understand a
continuous operator with norm bounded by a constant independent of
$\lambda$.

Now assume that $(A,B)$ acts in the half space $\R^n_+$ with
coefficients of the form \eqnref{2-1} and that $(A,B)$ is $N$-elliptic 
in the sense of Definition \ref{2.1}. To define a parametrix, we use
a cut-off function $\psi'\in C^\infty(\R^{n-1})$ with
\[ \psi'(\xi') = \left\{ \begin{array}{ll} 0\,,& |\xi'|\le 1\,,\\
1\,,& |\xi'| \ge 2\,.\end{array}\right.\]
For $j=1,\ldots,m$ we define the ps.d.o. $P_j$ in $\R^{n-1}$ (with
$x_n$ as parameter) by
\begin{equation}\mylabel{5-2}
(P_j g)(x',x_n) := \psi'(D')\, w_j(x',x_n,D',\lambda) \,g\,,
\end{equation}
where $w_j(x',x_n,\xi',\lambda)$ is the unique solution of
\eqnref{eq2-1}--\eqnref{eq2-2} with
\begin{equation}\mylabel{5-2a}
\begin{array}{rcl} A(\xi',D_t,\lambda) & = &
  A^{(0)}(x',0,\xi',D_t,\lambda)\,,\\[2mm]
 B_k(\xi',D_t) & = & B_k^{(0)}(x',0,\xi',D_t)\,.
\end{array}
\end{equation}
Due to Lemma \ref{2.3}, for
large $\lambda$ the symbol of $w_j(x',x_n,D',\lambda)$ can be written
in the form 
\begin{eqnarray}
w_j(x',x_n,\xi',\lambda) & = & \frac{1}{2\pi i}\int_{\gamma^{(1)}}
\frac{M_j^{(1)}(x',\xi',\tau,\lambda)}{A_1(x',\xi',\tau,\lambda)}\;
e^{i x_n\tau}\; d\tau \nonumber\\
& + & \frac{1}{2\pi i}\int_{\gamma^{(2)}}
\frac{M_j^{(2)}(x',\xi',\tau,\lambda)}{A_2(x',\xi'/\lambda,\tau,1)}\;
e^{i x_n\lambda\tau}\; d\tau \,.\mylabel{5-3}
\end{eqnarray}

\begin{lemma}\mylabel{5.2} The operator $P_j$ defined in \eqnref{5-2}
  is continuous from\\ $H^{\Xi^{(-m_j-1/2)}}(\R^{n-1})$ to
  $H^\Xi(\R^n_+)$.
\end{lemma}

\begin{proof} Let $g\in H^{\Xi^{(-m_j-1/2)}}(\R^{n-1})$ and set $u :=
  P_jg$. Using the equivalent norm
\begin{equation}\mylabel{5-3a}
 \left[ \sum_{l=0}^r \int_0^\infty \| D_n^l
  u(\cdot,x_n)\|^2_{\Xi^{(-l)},\R^{n-1}} \; dx_n\right]^{1/2} 
\end{equation}
in $H^\Xi(\R^n_+)$, we see that we have to show that
\[ \left\| \Xi^{(-l)}(D',\lambda) D_n^l P_j \left[
    \Xi^{(-m_j-1/2)}(D', \lambda)\right]^{-1}\right\|_{L_2(\R^{n-1})
  \to  L_2(\R^{n-1})} 
\le C(x_n)\]
for some function $C=C(x_n)$ whose $L_2(\R_+)$-norm is bounded by a
constant independent of $\lambda$. For this it is sufficient to show
that for $|\xi'|\ge 1$ we have
\[ \left( \int_0^\infty | D_{x'}^{\alpha'} D_n^l
  w_j(x',x_n,\xi',\lambda )|^2 \; dx_n\right)^{1/2} \le C\; 
\frac{\Xi^{(-m_j-1/2)}(\xi',\lambda)}{\Xi^{(-l)}(\xi',\lambda)}\;.\]
As we have for $|\xi'|\ge 1$ the equivalence
\[ \Xi_{r,s}(\xi',\lambda) \approx \Phi_{r,s}(\xi',\lambda)\]
for all $r,s\in\R$ (with $\Phi_{r,s}$ defined by the right-hand side
of \eqnref{3-4}), the case $\alpha'=0$ is already covered by Lemma
\ref{3.1}. Here we take into account that, due to condition
\eqnref{2-1}, the constant $C$ in Lemma \ref{3.1} applied to the
symbols \eqnref{5-2a} may be chosen
independently of $x'\in\R^{n-1}$. 

The case $\alpha'>0$ follows after differentiation of \eqnref{5-3}
with respect to $x'$ along the same lines as in the proof of Lemma
\ref{2.3}.
\end{proof}

\begin{lemma}\mylabel{5.3} The operator
\[ C_j := A(x,D,\lambda)P_j: H^{\Xi^{(-m_j-1/2)}}(\R^{n-1})\to
H^\Theta (\R^n_+)\]
is bounded. Here $\Theta(\xi,\lambda)$ is defined in \eqnref{5-1}.
\end{lemma}

\begin{proof}
The symbol of the ps.d.o. $C_j$ in $\R^{n-1}$ with parameter $x_n$ is
given by
\[ \sum_{|\alpha'|=1}^{2m}\frac 1{(\alpha')!}\partial^{\alpha'}_{\xi'}
A(x,\xi',D_n,\lambda)\;D^{\alpha'}_{x'}P_j(x',x_n,\xi',\lambda)\]
with
\[ P_j(x,\xi',\lambda)=\psi(\xi')\,w_j(x,\xi',\lambda)\,.\]
Consider the famliy ${\cal F} = \{ A(x,\xi,\lambda): x\in\R^n_+\}$ of
polynomials in $(\xi,\lambda)\in \R^{n+1}$. As the degree of the
polynomial $A(x,\cdot)$ is equal to $2m$ for all $x\in\R^n_+$, the
family $\cal F$ is a subset of the finite-dimensional vector space of
all polynomials in $(\xi,\lambda)$ of degree not greater than
$2m$. Therefore, there exists a finite set $x^{(1)},\ldots, x^{(K)}\in 
\R^n_+$ such that every $A\in{\cal F}$ may be represented in the form
\[ A(x,\xi,\lambda)=\sum_{k=1}^K c_k(x)\, A(x^{(k)},\xi,\lambda)\]
with smooth coefficients $c_k(x)$. 

Taking into account that the operators of multiplication by $c_k(x)$ are
bounded in $H^{\Theta}(\R^n_+)$, we reduce our problem to the proof of 
the boundedness of operators of the form
\[ C_{\alpha', l}:= a_{\alpha', l}(D',\lambda)\;
D_n^lD_{x'}^{\alpha'}P_j(x,D',\lambda):
H^{\Xi^{(-m_j-1/2)}}(\R^{n-1}) \rightarrow H^{\Theta}(\R^n_+)\,,\]
where
\begin{equation}\mylabel{5-4}
 |a_{\alpha', l}(\xi',\lambda)| \le C\;
 \Xi^{(-l-1)}_{2m,2\mu}(\xi',\lambda) \,.
\end{equation}     
Literally repeating the proof of Lemma \ref{5.2}
 we establish the boundedness of the
operator
\[ D^{\alpha'}_{x'}P_j(x,D,\lambda): 
H^{\Xi^{(-m_j-1/2)}}(\R^{n-1}) \rightarrow H^{\Xi}(\R^n_+)\,.\]
According to \eqnref{5-4} the operator 
\[ a_{\alpha',l}(D',\lambda): H^{\Xi}(\R^n_+) \rightarrow
H^{\Theta}(\R^n_+) \]
is bounded. As $C_{\alpha', l}$ is the product of of the above operators
this operator is also bounded.
\end{proof}

\begin{theorem}\mylabel{5.4} Consider in the half space $\R^n_+$ 
 the boundary value problem
  $(A,B)=(A,B_1,\ldots,B_m)$ of the form \eqnref{1-1}, \eqnref{1-4}
  with coefficients of the form
  \eqnref{2-1}. Assume that $A$ degenerates regularly at the
  boundary and that $(A,B)$ is
 $N$-elliptic with parameter in $[0,\infty)$ in the
  sense of Definition {\rm \ref{2.1}}. Then there exists a bounded operator
\[ P: H^{(r-2m,s-2\mu)}(\R^n_+) \times \prod_{j=1}^m
H^{\Xi^{(-m_j-1/2)}}(\R^{n-1}) \rightarrow H^{\Xi}(\R^n_+)\]
such that
\[ (A,B) P = I + T\]
where $I$ stands for the identity operator in the space
\begin{equation}\mylabel{5-8}
 H^{(r-2m,s-2\mu)}(\R^n_+) \times \prod_{j=1}^m
H^{\Xi^{(-m_j-1/2)}}(\R^{n-1})
\end{equation}
and $T$ is a continuous operator from the space \eqnref{5-8} to the space
\[ H^{\Theta}(\R^n_+) \times 
\prod_{j=1}^m H^{\Xi^{(-m_j+1/2)}}(\R^{n-1}) \]
with $\Theta(\xi,\lambda)$ being  defined in \eqnref{5-1}.
\end{theorem}

\begin{proof}
We define
\[ P(f,g_1,\ldots,g_m) := P_0 f + \sum_{j=1}^m P_j(g_j-B_jP_0f)\]
with $P_0$ from Lemma \ref{5.1} and $P_j$ given by \eqnref{5-2}. The
continuity of $P$ follows from Lemma \ref{5.1} and Lemma
\ref{5.2}. In order to see that the operator $T$ is continuous with
respect to the spaces given in the theorem, we denote the components
of $T$ by $T_0, T_1,\ldots, T_m$. The operator $T_0$ is given by
\[ T_0(f,g_1,\ldots,g_m) = AP_0 f - f +\sum_{j=1}^m AP_j(g_j-B_j
P_0f)\,.\]
We see from Lemma \ref{5.1} and Lemma \ref{5.3} that $T_0$ maps the
space \eqnref{5-8} continuously into $H^\Theta(\R^n_+)$. 

Turning to the other components $T_1,\ldots,T_m$, we remark that for
$j,k=1,\ldots,m$ the operator $B_k P_j$ equals $\delta_{kj}\,I$ 
up to operators of lower order. More precisely, the operator
\[ B_k(x,D) P_j - \delta_{kj} I \]
is a ps.d.o. in $\R^{n-1}$ which is continuous from
\[ H^{\Xi^{(-m_j-1/2)}}(\R^{n-1})\quad\mbox{ to }\quad
H^{\Xi^{(-m_k+1/2)}}(\R^{n-1})\,.\]
This is due to the fact that $w_j(x',x_n,\xi',\lambda)$ satisfies
\eqnref{eq2-1}--\eqnref{eq2-2}; the estimates for the lower order
terms of the ps.d.o. $B_kP_j$ can be found in the same way as it was
done for $AP_j$ in the proof of Lemma \ref{5.3}. From the continuity
of $B_k P_j - \delta_{kj} I$ the continuity of $T_k$ in the
spaces given in the theorem immediately follows.
\end{proof}

\mysection{6}{Proof of the necessity}

The aim of this section is to prove the following theorem.

\begin{theorem}\mylabel{4.1} Let $A$ degenerate regularly at the
  boundary  ${\partial M}$
  and assume that   inequality \eqnref{1-4a} holds. Let $r$ and $s$ be 
  integers satisfying \eqnref{3-2b} and assume, in addition, that
\begin{equation}\mylabel{4-0}
r\ge m \quad \mbox{ and }\quad \mu\le s \le r-m+\mu\,.
\end{equation}
If the a priori estimate \eqnref{3-6} holds, then $(A,B_1,\ldots,B_m)$ 
is $N$-elliptic with parameter in the sense of Definition {\rm
  \ref{2.1}}. 
\end{theorem}

Note that if the orders $m_j$ of the boundary operators $B_j$ are all
  different (e.g., if the boundary operators are normal), then
  \eqnref{1-4a} is  satisfied
and \eqnref{4-0} is a consequence of \eqnref{3-2b}.
The proof of this theorem is divided into several 
steps.

{\it Necessity of condition {\rm \ref{2.1} a)}}.
First of all note that applying estimate \eqnref{3-6}
 to functions whose 
support does not  intersect with the boundary, we obtain the estimate in
$\R^n$
\begin{equation}\mylabel{ne1-1}
||u||_{r,s} \le C\Big(||A(x,D,\lambda)u||_{r-2m,s-2\mu}
+\lambda^{r-s}||u||_{L_2}\Big) \,     
\end{equation}

\begin{proposition}\mylabel{ne1.1}
Suppose that \eqnref{ne1-1} takes place and $x^0$ is an arbitrary 
point of $\R^n$. Denote $A(D,\lambda)=A^0(x^0,D,\lambda)$. Then the estimate
\begin{eqnarray}
\lefteqn{\Big\| |D|^{s}(|D|+\lambda)^{r-s}u \Big\|_{L_2}}
\nonumber\\
 &\le& C\Big\| |D|^{s-
2\mu}(|D|+\lambda)^{r-s-2m+2\mu}A(D,\lambda)u\Big\|_{L_2} \mylabel{ne1-2}
\end{eqnarray}
holds, where the constant $C$ does not depend on $x^0$ or $\lambda$.
\end{proposition}

The necessity of a) easily follows from \eqnref{ne1-2}. Indeed,
applying the Fourier transform,
we can rewrite \eqnref{ne1-2} in the form
\begin{eqnarray*}
\lefteqn{
\int_{\R^n}\Big[|\xi|^{2s}
(|\xi|+\lambda)^{2(r-s)}}\\
& & - C^2
|\xi|^{2(s-2\mu)}
 (|\xi|+\lambda)^{2(r-s-2m+2\mu)}|A(\xi,\lambda)|^2\Big]
|(Fu)(\xi)|^2 d\xi \le 0\,.
\end{eqnarray*}
Since $u\in \Cal D$ is arbitrary, the expression in the square brackets is
nonpositive. From this part  a) follows.

To prove \eqnref{ne1-2} we replace in \eqnref{ne1-1}
 $\lambda$  by $\rho \lambda$ with $\rho>0$ and $u(x)$ by 
\begin{equation}\mylabel{ne1-3}
u_{\rho}(x)=\rho^{-r+n/2}u(\rho(x-x^0))     
\end{equation}
and tend $\rho$ to $+\infty$. To carry out the calculations we need the
following

\begin{lemma}\mylabel{ne1.2}
Denote
\begin{equation}\mylabel{ne1-4}
(S_{\rho,x^0}u)(x)=u(\rho(x-x^0)).  
\end{equation}
Then for an arbitrary ps.d.o. $a(x,D)$ we have
\begin{equation}\mylabel{ne1-5}
\Big[a(x,D)S_{\rho,x^0}u\Big](x)=\Big[ 
S_{\rho,x^0} a(x^0+\rho^{-1}x,\rho D)u\Big](x)\,.
  \end{equation}                                                      
\end{lemma}

\begin{proof}
Direct calculation shows that
\[
(F S_{\rho,x^0}u)(\xi)=\rho^{-n}\exp(-ix^0\xi)(Fu)\Big(
\frac {\xi}{\rho}\Big)\,.
\]
If we substitute the last expression in the left-hand side
 of \eqnref{ne1-5} and 
change $\xi$ to $\rho\xi$ we obtain the right-hand side of
\eqnref{ne1-5}.
\end{proof}

\begin{proofof}{Proof of Proposition {\rm \ref{ne1.1}}} Applying the a 
  priori estimate \eqnref{ne1-1} to the function $u_\rho$
  (cf. \eqnref{ne1-3}), we obtain, according to the lemma,
\begin{eqnarray*}
\lefteqn{((1+|D|)^{s}(\rho\lambda+|D|)^{r-s}u_{\rho})(x)}\\
&&=\rho^{n/2}
S_{\rho,x^0}\Big[(\rho^{-1}+|D|)^{s}(\lambda + |D|)^{r-s}u
\Big](x)\,.
\end{eqnarray*}
The $L_2(\R^n)$ norm of this expression tends to the left-hand 
side of \eqnref{ne1-2},
as $\rho$ tends to $+\infty$.

Now we turn to the right-hand side of \eqnref{ne1-1}. We have
\begin{eqnarray*}
\lefteqn{ 
(1+|D|)^{s-2\mu}(\rho\lambda+|D|)^{r-s-2m+2\mu}A(x,D,\rho\lambda)
u_{\rho}(x) } \\
&&  = \rho^{n/2}S_{\rho,x^0}\Big[(\rho^{-1}+|D|)^{s-2\mu}
(\lambda+|D|)^{r-s-2m+2\mu}h_{\rho}\Big](x)
\end{eqnarray*}
where
\[
h_{\rho}(x)=\rho^{-2m}A(x^0+\rho^{-1}x,\rho D,\rho \lambda)u(x) 
\quad\rightarrow \quad A(D,\lambda)u
\]
as $\rho \rightarrow +\infty$.

It is easy to check that the limit of the second term of the right-hand side 
of \eqnref{ne1-1} is equal to zero. 
\end{proofof}  

\vskip 5 mm

To prove the necessity of \ref{2.1} b), 
c) and  d) we consider \eqnref{3-6} for functions
with supports belonging to a small neighbourhood of a point $x^0\in
{\partial M}$. In this case the norms in \eqnref{3-6} can be taken in $\R^n_+$
and  $\R^{n-1}$, respectively. Now we use the fact that we have the
norm equivalence \eqnref{5-3a}
 and the equivalence
\[ \Xi^{(-l)}(\xi,\lambda) \approx \left\{ \begin{array}{ll}
(1+|\xi|)^{s-l} \, (\lambda+|\xi|)^{r-s}\,,& 
l\le s\,,\\
(\lambda+|\xi|)^{r-l}\,,& l> s\,.
\end{array}\right. \]
According to  \cite{dmv2}, Section 2, 
 the norm 
\[ \| (iD_n + \sqrt{ 1+|D'|^2})^q \, (iD_n + \sqrt{\lambda^2 +
  |D'|^2}) ^{p-q} u \|_{L_2(\R^n_+)} \]
is defined for any $p,q\in\R$ and 
is equivalent to $\|u \|_{(p,q),\R^n_+}$
(cf. \eqnref{3-2a}). Substituting these expressions into the a priori
estimate \eqnref{3-6} for the half space, we obtain in explicit form
\begin{eqnarray} 
\lefteqn{\sum_{l=0}^{s}\Big\|
(1+|D'|^2)^{(s-l)/2}(\lambda^2+
|D'|^2)^{(r-s)/2}D_n^lu\Big\|_{L_2(\R^n_+)}}\nonumber\\
\lefteqn{+ \sum_{l=s+1}^{r}\Big\|(\lambda^2+|D'|^2)^{(r-l)/2}
D_n^lu\Big\|_{L_2(\R^n_+)} }\nonumber\\
& \le & C \Bigg(\|\sigma(D,\lambda)A(x,D,\lambda)u\|_
{L_2(\R^n_+)}\nonumber\\
& + & \lambda^{r-s}||u||_{L_2(\R^n_+)}\nonumber\\
&  + & \!\!\sum_{j=1}^{\mu}\Big\|
(1+|D'|^2)^{(s-m_j-1/2)/ 2}
(\lambda^2+|D'|^2)^{(r-s)/2}B_j(x',D)u\Big\|_{L_2(\R^{n-1})}
\nonumber\\
&  + &\!\! \sum_{j=\mu+1}^{m}\Big\|(\lambda^2+|D'|^2)^{(r-m_j-1/2)/2}
B_j(x',D)u\Big\|_{L_2(\R^{n-1})}\Bigg) \,,   \mylabel{ne2-1}           
\end{eqnarray}
where we used the abbreviation
\[ \sigma(D,\lambda) := (iD_n+\sqrt{1+|D'|^2})^{s-2\mu}
(iD_n+\sqrt{\lambda^2+|D'|^2})^{r-s-2m+2\mu}\,.\]

\begin{proposition}\mylabel{ne2.2}
 Suppose estimate \eqnref{ne2-1}
 holds. Let $x^0$ be an arbitrary
point in $\R^{n-1}$ and set $A(D,\lambda)=A^{(0)}(x^0,D,\lambda),\quad
B_j(D,\lambda)=B_j^{(0)}(x^0,D,\lambda),\; j=1,\dots,m$. Then the
following estimate holds
\begin{eqnarray}
\lefteqn{\sum_{l=0}^{s}\Big\||D'|^{s-l}(\lambda^2+|D'|^2)^
{(r-s)/2}D_n^lu\Big\|_{L_2(\R^n_+)} }\nonumber\\
\lefteqn{+\sum_{l=s+1}
^{r}\Big\|(\lambda^2+|D'|^2)^{(r-l)/2}
D_n^lu\Big\|_{L_2(\R^n_+)}}\nonumber\\
& \le & C\Bigg(\| \tilde\sigma(D,\lambda) A(D,\lambda)u\|_
{L_2(\R^n_+)}\nonumber\\
&+&\sum_{j=1}^{\mu}\Big\| |D'|^{s-m_j-1/2}
(\lambda^2+|D'|^2)^{(r-s)/2}B_j(D)u\Big\|_{L_2(\R^{n-1})}
\nonumber\\
& + & \sum_{j=\mu+1}^{m}\Big\|(\lambda^2+|D'|^2)^{(r-m_j-1/2)/2}
B_j(D)u\Big\|_{L_2(\R^{n-1})}\Bigg)\mylabel{ne2-2} \,,           
\end{eqnarray}
where we have set
\[ \tilde\sigma(D,\lambda) := (iD_n+|D'|)^{s-2\mu}
(iD_n+\sqrt{\lambda^2+|D'|^2})^{r-s-2m+2\mu}\,.\]

\end{proposition}

\begin{proof} We apply \eqnref{ne2-1} with $\lambda$ replaced by
  $\rho\lambda$ to the function $u_\rho$ defined in \eqnref{ne1-3},
  noting that $S_{\rho,x^0}u$ is again defined in $\R^n_+$ because of
  $x^0\in\R^{n-1}$ and $\rho>0$. From Lemma \ref{ne1.2} and the fact
  that for any function $v\in L_2(\R^n_+)$ we have
\[ \rho^{n/2} \| S_{\rho, x^0} v\|_{L_2(\R^n_+)} = 
\| v \|_{L_2(\R^n_+)}\,,\]
we see that the $l$-th term in the first sum in \eqnref{ne2-1} is
equal to 
\[ \Big\| \Big( \frac 1{\rho^2} + |D'|^2\Big)^{(s-l)/2}\,
\Big( \lambda^2 + |D'|^2\Big)^{(r-s)/2} D_n^l u\Big\|_{L_2
  (\R^n_+)} \]
which tends to the corresponding term in \eqnref{ne2-2} for $\rho \to
\infty$. The remaining expressions in \eqnref{ne2-1} can be treated
analogously; the term $(\rho\lambda)^{r-s}
\|u_\rho\|_{L_2(\R^n_+)}$ tends to zero for $\rho\to\infty$.

For the terms involving the boundary operators we remark that
$\gamma_0 S_{\rho,x^0} = S_{\rho, x^0}\gamma_0$ where $\gamma_0:
u\mapsto u(\cdot,0)$ stands for the trace operator. Therefore we may
apply Lemma \ref{ne1.2} to the function $B_j(x',D)u_\rho$ defined in
$\R^{n-1}$. 
\end{proof}

If we apply \eqnref{ne2-2} to a function of the form
\[
u(x)=\phi(x')V(x_n),\quad \phi(x') \in \Cal D(\R^{n-1})
\]
we obtain an estimate on the half-line (cf.
\cite{gindikin-volevich2}, Chapter 3, Proposition 2 in Subsection 2.3):

\begin{eqnarray}
\lefteqn{
\sum_{l=0}^{s}|\xi'|^{s-l}(\lambda^2+|\xi'|^2)^
{(r-s)/2}||D_n^lV||_{L_2(\R_+)}}\nonumber\\
&+&\sum_{l=s+1}
^{r}(\lambda^2+|\xi'|^2)^{(r-l)/2}||D_n^lV||_{L_2(\R_+)} \nonumber\\
&\le& C
\Bigg(\| \tilde\sigma(\xi',D_n,\lambda)
 A(\xi',D_n,\lambda)V\|_{L_2(\R_+)} \nonumber\\
&+&\sum_{j=1}^{\mu}
|\xi'|^{s-m_j-1/2}(\lambda^2+|\xi'|^2)^{(r-s)/2}
|B_j(\xi',D_n)V(0)| \nonumber\\
&+& \sum_{j=\mu+1}^{m}(\lambda^2+|\xi'|^2)^{(r-m_j-1/2)/2}
|B_j(\xi',D_n)V(0)|\Bigg) \,. \mylabel{ne2-3}           
\end{eqnarray}

\vskip 5mm

{\it     Necessity of condition {\rm \ref{2.1} b)}}.
Suppose $V(x_n)\in L_2(\R_+)$ is a solution of the homogeneous equation
\[
A(\xi',D_n,\lambda)V(x_n)=0,\quad x_n>0\,.
\]
Then this function satisfies the equation
\begin{equation}\mylabel{ne2-4}
A_+(\xi',D_n,\lambda)V(x_n)=0,\quad x_n>0. 
\end{equation}
Now from \eqnref{ne2-3} we deduce the estimate
\begin{equation}\mylabel{ne2-5}
c(\xi',\lambda)\sum_{l=0}^{r}||D_n^l V||_{L_2(\R_+)} \le
\sum_{j=1}^m |B'_j(\xi',\lambda,D_n)V(0)|                       
\end{equation}
Here $B'_j$ are  remainders of $B_j$ after the division by $A_+$ and
$c(\xi,\lambda)>0$ for $\xi' \ne 0$ and $\lambda \ge 0$. From a
standard trace result for Sobolev spaces on $\R_+$ we know that
\begin{equation}\mylabel{4-10a}
  \sum_{j=1}^r | D_n^{j-1} V(0)| \le C \sum_{j=1}^{r+1} \| D_n^{j-1}
  V\|_{L_2(\R_+)}\,.
\end{equation}
\mbox{}From this and  \eqnref{ne2-5} we obtain, using $r\ge m$ (see
\eqnref{4-0}), 
\begin{equation}\mylabel{4-10b}
\tilde c(\xi',\lambda)\sum_{j=1}^m|D_n^{j-1}V(0)| \le \sum_{j=1}^m
\Big|\sum_{k=1}^m b_{jk}(\xi',\lambda)D_n^{k-1}V(0)\Big|,
\end{equation}
where 
\begin{equation}\mylabel{ne2-6}
B'_j(\xi',\lambda,z)=\sum_{k=1}^m b_{jk}(\xi',\lambda)z^{k-1}\,.
\end{equation}       
The constant $\tilde c(\xi',\lambda)$ in \eqnref{4-10b} is positive
for $\xi'\not=0$ and $\lambda\ge 0$.       
Note that the Cauchy problem
\[
D_n^{k-1}V(0)=\zeta_k,\quad k=1,\dots,m
\]
for ODE \eqnref{ne2-4}
 has a unique solution for arbitrary $\zeta=(\zeta_1,\dots,
\zeta_m) \in \C^m$. This means that for an arbitrary complex vector
$\zeta$ we have the estimate
\[
\tilde c(\xi',\lambda)|\zeta| \le |\Cal B(\xi',\lambda)\zeta|
\]
where $\Cal B(\xi',\lambda):=(b_{jk}(\xi',\lambda))_{j,k=1,\dots,m}$.
The last inequality means that the matrix $\Cal B(\xi',\lambda)$ is 
nonsingular as $|\xi'| \ne 0,\lambda \ge 0$, i.e. the necessity of b)
is proved.

\begin{proposition}\mylabel{ne2.3}
 Suppose the estimate \eqnref{ne2-1} holds. Let $x^0$
be an arbitrary point of $\R^{n-1}$. Then the inequality 
\begin{eqnarray}
\lefteqn{
\sum_{l=0}^{s}\Big\| |D'|^{s-l}D_n^lu\Big\|_{L_2(\R^n_+)}}
\nonumber\\
& \le &
C\Bigg(\Big\|(iD_n+|D'|)^{s-2\mu}A_{2\mu}(D)u\Big\|_{L_2(\R^n_+)}
\nonumber\\
&+&\sum
_{j=1}^{\mu}\Big\| |D'|^{s-m_j-1/2}B_j(D)u\Big\|_{L_2(\R^{n-1})}
\Bigg)
\mylabel{ne2-6a}
\end{eqnarray}
holds, where  $A_{2\mu}(D)=A^{(0)}_{2\mu}(x^0,D),\; B_j(D)=
B^{(0)}_j(x^0,D),\,j=1,\ldots,\mu$.
\end{proposition}

\begin{proof} This can be seen in exactly the same way as Proposition
  \ref{ne2.2}, now applying the a priori estimate \eqnref{ne2-1} with
  $\rho^t\lambda$ instead of $\lambda$ to the function
\[ u_\rho(x) := \rho^{-t(r-s)-s+n/2}
u(\rho(x-x^0)) \] 
where $t>1$ is fixed and $\rho>0$ tends to infinity.
\end{proof}

{\it Necessity of condition {\rm \ref{2.1} c)}}.
From Proposition \ref{ne2.3} 
 the estimate on the half-line can be obtained
\begin{eqnarray*}
\lefteqn{
\sum_{l=0}^{s}|\xi'|^{s-l}||D_n^lV||_{L_2(\R_+)}}
\nonumber\\
& \le &
C\Bigg(\Big\|(iD_n+|\xi'|)^{s-2\mu}A_{2\mu}(\xi',D_n) 
V\Big\|_{L_2(\R_+)}\nonumber\\
&+&\sum
_{j=1}^{\mu}|\xi'|^{s-m_j-1/2}|B_j(\xi',D_n)V(0)|\Bigg)\,.
\end{eqnarray*}
As above we see that for solutions $V(x_n)\in L_2(\R_+)$ of
\[
A_{2\mu}(\xi',D_n)V(x_n)=0, \quad x_n>0
\]
we obtain the inequality
\begin{equation}\mylabel{ne2-7}
\sum_{l=1}^{s} |D_n^{l-1}V(0)| \le C\sum_{j=1}^{\mu}|B'_j
(\xi',\lambda,D_n)V(0)|    
\end{equation}
with a constant $C$ independent of $\xi',\,|\xi'|=1$, and $\lambda$,
where now $B_j'$ denotes the remainder of $B_j$ after division by
$(A_{2\mu})_+$. 
Replacing in \eqnref{ne2-7}
 the germ of $V$ in $0$ by an arbitrary vector $\zeta \in
\C^{\mu}$ and using $s\ge \mu$ (see \eqnref{4-0}),
 we obtain the necessity of c). 

\begin{proposition}\mylabel{ne2.4} 
Suppose the estimate \eqnref{ne2-1} holds and $x^0$ is an
arbitrary point of $\R^{n-1}$. Then the estimate
\begin{eqnarray}
\sum_{l=s}^{r} ||D_n^lu||_{L_2(\R^n_+)}
& \le & C\Bigg(||(D_n-i)^{r-s-2m+2\mu}
D_n^{s}Q(x^0,D_n)u||_{L_2(\R^n_+)} \nonumber\\
&+&\sum_{j=\mu+1}^m||B_j^{(0)}(x^0,0,D_n)u||
_{L_2(\R^{n-1})}\Bigg)\mylabel{4-14}
\end{eqnarray}
holds.
\end{proposition}

\begin{proof} We apply \eqnref{ne2-1} with $\lambda$ replaced by
  $\rho$ to the function
\[ u_\rho(x) := \rho^{1/2 + \epsilon (n-1)/2 - r} u\big(
\rho^\epsilon(x'-x^0), \rho x_n\big) \]
with $0<\epsilon <1$ fixed. Now we use
\[ u\big( \rho^\epsilon(x'-x^0), \rho x_n\big) = \Big[
S_{\rho^\epsilon, x^0}^{(x')} \; S_{\rho,0}^{(x_n)} u \Big] (x) \,,\]
where $S_{\cdots}^{(x')}$
 indicates that the operator $S_{\cdots}$ acts
on the first $n-1$ variables (and analogously that $S_{\cdots}^{(x_n)}$
acts on the last variable), and apply Lemma \ref{ne1.2} twice. For
the $l$-th term in the first sum of \eqnref{ne2-1} we obtain the
expression
\[ \rho^{(1-\epsilon)(l-s)} \Big\| (\rho^{-2\epsilon} +
|D'|^2)^{(s-l)/2}  (1 + \rho^{2(\epsilon-1)}
|D'|^2)^{(r-s)/2}  D_n^l u \Big\|_{L_2(\R^n_+)}\!. \]
For $l\le s-1$ this expression tends to zero for
$\rho\to\infty$, for $l=s$ its limit equals $\|
D_n^{s} u\|_{L_2(\R^n_+)}$. 

The remaining terms can be treated analogously; to finish the proof we 
use
\[ \rho^{-2m} \, A\Big( x^0 + \frac{x'}{\rho^\epsilon},
\frac{x_n}{\rho}, \rho^\epsilon D', \rho D_n, \rho\Big) \; \to \;
D_n^{2\mu} Q(x^0,D_n)\]
and
\[ \rho^{-m_j} B_j\Big(  x^0 + \frac{x'}{\rho^\epsilon},
 \rho^\epsilon D', \rho D_n \Big)\; \to \;
B_j^{(0)}(x^0,0,D_n) \]
as $\rho\to\infty$.
\end{proof}

{\it Necessity of condition {\rm \ref{2.1} d)}}.
From Proposition \ref{ne2.4} we obtain  the estimate on the half-line
\begin{eqnarray}
\sum_{l=s}^{r} ||D_n^lV||_{L_2(\R_+)} 
&\le& C\Bigg(\Big\|(D_n-i)^{r-s-2m+2\mu}
D_n^{s}Q(x^0,D_n)V\Big\|_{L_2(\R_+)} \nonumber\\
&+&\sum_{j=\mu+1}^m|B_j^{(0)}(x^0,0,D_n)V(0)|\Bigg)\,.\mylabel{4-15}
\end{eqnarray}
Since $m_j\ge s$ for $j\ge\mu+1$, 
each $B_j^{(0)}(x^0,0,\tau)$ contains the factor $\tau^{s}$,
and it is easily seen that condition d) follows from the analogous
condition for 
\[ \tilde B_j^{(0)}(x^0,0,\tau) := \tau^{-s}
B_j^{(0)}(x^0,0,\tau) \,.\]
Now we apply \eqnref{4-15} to a  solution $V\in L_2(\R_+)$ of 
\[
Q(x^0,D_n)V(x_n)=0,\quad x_n>0
\]
and substitute $W(x_n) := D_n^{s} V(x_n)$. We obtain
\begin{equation}\mylabel{4-16}
\sum_{l=1}^{r-s+1}\| D_n^{l-1} W\|_{L_2(\R_+)}
 \le C\sum_{j=\mu+1}^m|B'_j(x^0,0,D_n) W(0)|,
\end{equation}
where now $B'_j$ stands for the 
 remainder of $\tilde B_j^{(0)}$ after division  by $Q_+$. Using
 $r-s  \ge m-\mu$ (cf. \eqnref{4-0}) and the trace result 
 \eqnref{4-10a}, we obtain the linear independence of $B_j'$ modulo
 $Q_+$ from \eqnref{4-16} and therefore condition d).

\myappendix{Singularly perturbed problems}

One of the most important features of the Newton polygon approach is
to provide an easy formulation and proof of a priori estimates in the
theory of singularly perturbed problems. All results of the previous
sections can be rewritten for boundary value problems with small
parameter as treated by Vishik--Lyusternik \cite{vishik-lyusternik},
Nazarov \cite{nazarov}, Frank \cite{frank} and others. (Cf. also
\cite{dmv2}, Appendix, for the Dirichlet problem.) As an example, we
formulate an a priori estimate for such problems.

Consider for $\epsilon>0$ the operator
\[
A_\epsilon(x,D) := \epsilon^{2m-2\mu} A_{2m}(x,D) +
\epsilon^{2m-2\mu-1} A_{2m-1}(x,D) + \ldots + A_{2\mu}(x,D)
\]
with $A_j$ of the form \eqnref{1-2}. Let $A_\epsilon$ act on a smooth
compact manifold $M$ with boundary ${\partial M}$ and assume that we have
boundary conditions $B_1(x,D),\ldots,B_m(x,D)$ of the form
\eqnref{1-4} satisfying \eqnref{1-4a}. 

We fix integer numbers $r$ and $s$ satisfying \eqnref{3-2b} and
consider the weight function
\[ \Xi_\epsilon(\xi) := \Xi_{\epsilon, (r,s)}(\xi) := (1+|\xi|)^s \;
(1+\epsilon |\xi|)^{r-s}\,.\]
The norms corresponding to this weight function will be denoted by\\
$\|\cdot \|_{\Xi_\epsilon,M} = \|\cdot\|_{\epsilon,(r,s),M}$.

\begin{definition}\mylabel{A.1}{\rm a) The operator $A_\epsilon(x,D)$
    is called $N$-elliptic if 
\[  |A_\epsilon^{(0)}(x,\xi)| \ge C |\xi|^{2\mu}\;(1+\epsilon
|\xi|)^{2m-2\mu}\quad (\xi\in\R^n\,,\;\epsilon>0\,,\; x\in\overline
M)
\]
holds where $C$ does not depend on $x,\xi$ or $\epsilon$.

\vskip 0.5em
\noindent b) The operator $A_\epsilon$ is said to 
degenerate regularly at the boundary if the
polynomial
\[ Q(x^0,\tau) := \tau^{-2\mu} A_1^{(0)}(x^0,0,\tau)\]
has exactly $m-\mu$ roots in the upper half plane. }
\end{definition}

\begin{definition}\mylabel{A.2}{\rm  The boundary problem
    $(A_\epsilon,B_1,\ldots, B_m)$ is called $N$-elliptic 
     if the following conditions hold:

\vskip 0.5em
\noindent a) The operator $A_\epsilon(x,D)$ is $N$-elliptic 
 in the sense of Definition \ref{A.1}.

\vskip 0.5em
\noindent b) For every fixed $x^0\in{\partial M}$ the boundary problem
\[ \Big( A^{(0)}_{\epsilon}(x^0,\xi',D_n),\;
B^{(0)}_1(x^0,\xi',D_n),\ldots,B^{(0)}_m(x^0,\xi',D_n)\,\Big)\]
for each $\epsilon>0$ and $\xi'\ne 0$ is uniquelly solvable on the half-line
$x_n\ge 0$ in the space of functions tending to zero as $x_n\rightarrow
\infty$. Moreover we suppose that the problem
\[ \Big(A^{(0)}_{2m}(x^0,\xi',D_n),\;
B^{(0)}_1(x^0,\xi',D_n),\ldots,B^{(0)}_m(x^0,\xi',D_n)\,\Big)\]
(corresponding to $\epsilon=\infty$) has the same property.

\vskip 0.5em
\noindent c) For every  $x^0 \in {\partial M}$ the
 boundary problem 
\[ (A_{2\mu}(x^0,D), B_1(x^0,D),\ldots,
  B_{\mu} (x^0,D))\]
 fulfills the Shapiro--Lopatinskii condition.

\vskip 0.5em
\noindent d) For every $x^0\in{\partial M}$ the  polynomials
  $(B_j^{(0)}(x^0,0,\tau))_{j=\mu+1,\ldots, m}$ are linearly independent 
  modulo $Q_+(x^0,\tau)$ with $Q_+$ defined in Definition \ref{2.1} d).

}
\end{definition}

If the conditions of Definition \ref{A.1} and \ref{A.2} hold, we can
apply Theorem \ref{3.2} to the operator
\[ A(x,D,\lambda) := \lambda^{2m-2\mu} \, A_{1/\lambda}(x,D)\,.\]
The connection between $\Xi_\epsilon(\xi)$ and $\Xi(\xi,
\epsilon^{-1})$ (defined in \eqnref{3-2}) is given by
\[ \Xi_\epsilon(\xi) = \epsilon^{r-s} \; \Xi(\xi,\epsilon^{-1})\]
and
\[ \Xi_\epsilon^{(-m_j-1/2)}(\xi) = \left\{ \begin{array}{ll}
\epsilon^{r-s} \, \Xi^{(-m_j-1/2)}(\xi,\epsilon^{-1})&\mbox{ if }j\le
\mu\,,\\ 
\epsilon^{r-m_j-1/2}\, \Xi^{(-m_j-1/2)}(\xi,\epsilon^{-1}) &
\mbox{ if } j>\mu\,. \end{array}\right. \]
Using these relations, we obtain from Theorem \ref{3.2} the following
result which can be found (without the notation of the Newton polygon) 
in \cite{frank}:

\begin{theorem}\mylabel{A.3} Assume that $A_\epsilon$ degenerates
  regularly and that\\ $(A_\epsilon, B_1,\ldots,B_m)$ is $N$-elliptic in 
  the sense of Definition {\rm \ref{A.2}}. Then the following 
a priori estimate holds with a constant $C$ independent of  $\epsilon > 0$:
\begin{eqnarray*}
\|u\|_{\Xi_\epsilon,M} & \le & C \Bigg( \| A_\epsilon u\|_{\epsilon,
(r-2m, s-2\mu), M} + \sum_{j=1}^\mu 
 \| B_j u\|_{\Xi_\epsilon^{(-m_j-1/2)},{\partial M}}\nonumber\\
& +&\sum_{j=\mu+1}^m \epsilon^{m_j+1/2-s}\;
 \| B_j u
\|_{\Xi_\epsilon^{(-m_j-1/2)},{\partial M}} +
  \|u\|_{L_2(M)}\Bigg)\,.
\end{eqnarray*}
\end{theorem}

\def\refname

\end{document}